

\documentstyle{amsppt}
\magnification=\magstep1
\NoRunningHeads

\vsize=7.4in


\def\supp{\text{supp }}


\topmatter

\title The basic sequence problem
\endtitle
\author
N.J. Kalton
\endauthor
\address
Department of Mathematics,
University of Missouri,
Columbia, MO  65211, U.S.A.
\endaddress
\email mathnjk\@mizzou1.missouri.edu \endemail
\thanks The author was supported in part by NSF-grant DMS-9201357
\endthanks
\abstract
We construct a quasi-Banach space $X$ which contains no basic sequence.
\endabstract

\endtopmatter

\document
\baselineskip=14pt

\heading{1. Introduction}\endheading \vskip10pt

It is a classical result in Banach space theory, known to Banach himself
\cite{1}, that every (infinite-dimensional) Banach space contains a
closed linear subspace with a basis, or, in other words, a basic
sequence.  The corresponding question for quasi-Banach spaces (and more
general F-spaces) has, however, remained open.  A number of equivalent
formulations are known (\cite{11},\cite{14},\cite{16},\cite{17}); the
question is also raised in a slightly disguised form in \cite{28} p.114.

In \cite{11} and \cite{17} it is shown that a quasi-Banach space $X$
contains a basic sequence if and only if there is a strictly weaker
Hausdorff vector topology on $X$.  Thus the existence of a space with no
basic sequence is equivalent to the existence of a (topologically) {\it
minimal} space (i.e. one on which there is no strictly weaker Hausdorff
vector topology).  See \cite{3} and \cite{4} for a discussion of minimal
spaces.  It further follows that $X$ contains a basic sequence if and
only if there is some infinite-dimensional closed subspace with
separating dual (\cite{11} Theorem 4.4).  Several positive results are
known.  For example, the work of Bastero \cite{2} implies that every
subspace of $L_p[0,1]$ ($0<p<1$) contains a basic sequence, while the
author's results in \cite{12} imply that every quotient of $L_p[0,1]$
contains a basic sequence.  Bastero's result can be lifted to the wider
class of so-called natural spaces and has further been extended by Tam
\cite{30} who shows that every complex quasi-Banach space with an
equivalent plurisubharmonic norm contains a basic sequence.  These
results suggest that almost all ``reasonable'' spaces contain a basic
sequence.

In this paper, we will prove:

\proclaim{Theorem 1.1} There is a quasi-Banach space $Y$ with a
one-dimensional subspace $L$ so that:\newline (1) If $Y_0$ is a closed
infinite-dimensional subspace of $Y$ then $L\subset Y_0.$\newline (2)
$Y/L$ is isomorphic to the Banach space $\ell_1.$ \newline In particular
$Y$ contains no basic sequence and is minimal.  \endproclaim

It is clear that (1) would make it impossible for $Y$ to contain a basic
sequence.

There are other applications of this space.  A topological vector space
$X$ is said to have the Hahn-Banach Extension Property (HBEP) if whenever
$X_0$ is a closed subspace of $X$ and $f$ is a continuous linear
functional on $X_0$ then $f$ can be extended to a continuous linear
functional on $X$.  The author showed in \cite{11}, answering a question
raised by Duren, Romberg and Shields \cite{5} (see also \cite{25},
\cite{29}) that for an F-space (complete metric linear space) (HBEP) is
equivalent to local convexity.  It was very well-known that metrizability
is necessary in this theorem, but some partial results of Ribe \cite{25}
suggested that completeness might not be required.  Ribe showed that if
$X$ is a metric linear space so that $X$ is isomorphic to $X\oplus X$
then if $X$ has (HBEP) it must be locally convex.  More recently, the
author \cite{14} extended Ribe's result to show:

\proclaim{Theorem 1.2}Let $X$ be a decomposable quasi-Banach space (i.e.
there is bounded projection $P$ on $X$ so that neither $P$ nor $I-P$ has
finite rank).  Suppose $X_0$ is a dense subspace of $X$.  Then $X_0$ has
(HBEP) if and only if $X$ is locally convex.\endproclaim

A proof of Theorem 1.2 is included in Section 6. The Hahn-Banach
extension property for metrizable spaces is also discussed in \cite{10}.

However, if $Y$ is the space constructed above, we will show that any
algebraic complement $Y_0$ of $L$ has (HBEP).  Thus we have:

\proclaim{Theorem 1.3}There is a non-locally convex metric linear space
$Y_0$ with the Hahn-Banach Extension Property.\endproclaim

In 1962, Klee \cite{18} asked whether for every topological vector space
$(X,\tau),$ the topology $\tau$ can be expressed as the supremum of two
not necessarily Hausdorff vector topologies $\tau_1$ and $\tau_2$ so that
(the Hausdorff quotient of) $(X,\tau_1)$ has a separating dual (i.e. is
{\it nearly convex}) and $(X,\tau_2)$ has trivial dual.  Recently Peck
\cite{22} has shown this to be true for certain twisted sums of a Banach
space and a one-dimensional space (see also \cite{23}).  The space
constructed here, $Y$, turns out to be a counterexample to Klee's
problem.

\proclaim{Theorem 1.4}There is a quasi-Banach space $Y$ so that the
topology on $Y$ is not the supremum of a trivial dual topology and a
nearly convex topology.\endproclaim

The construction of our example depends heavily on the recent remarkable
developments in infinite-dimensional Banach spaces due to Gowers, Maurey,
Odell and Schlumprecht \cite{7},\cite{8},\cite{9},\cite{20},\cite{21}.
It is perhaps a little ironic that the basic sequence question for
quasi-Banach spaces turns out to be so closely related to the {\it
unconditional} basic sequence problem for Banach spaces.  However, it
should be stressed that we use an example of a Banach space with an
unconditional basis, very similar to that used by Gowers in \cite{7}; the
fundamental estimates we need are in \cite{9}.

Let us conclude this introduction by explaining the shortcomings of the
example.  It is still an open question whether every quasi-Banach space
(or F-space) must contain a proper closed infinite-dimensional subspace.
A space with no proper closed infinite-dimensional subspace is called
{\it atomic}.  The existence of an atomic quasi-Banach space is known to
be equivalent to the existence of a {\it quotient minimal} quasi-Banach
space, i.e. a space $X$ so that every quotient is minimal (this concept
is due to Drewnowski \cite{3}).  See \cite{14} or \cite{16} for a
discussion.  Our example is quite far from an atomic space, and it is not
clear at the present whether it can be used towards making such a
monster.  We remark that Reese \cite{24} has constructed an example of an
``almost'' atomic F-space, i.e. a space $X$ with a sequence of
finite-dimensional subspaces $V_n$ with dim $V_n>n$ so that if $x_n\in
V_n$ is any sequence which is nonzero infinitely often then $[x_n]=X.$ It
is still unknown whether even this phenomenon can be reproduced in a
quasi-Banach space.  We suspect, however, that an atomic quasi-Banach
space will eventually be found.

We would like to thank several colleagues for helpful comments and
remarks during the course of this work, in particular P. Casazza, D.
Kutzarova, M. Lammers, M. Masty\l o and N.T.  Peck.  We also want to
thank B. Maurey for a substantial simplification of the last part of the
argument which we have incorporated into the proof.  We also wish to
thank the referee for many very helpful suggestions and comments on
improving the presentation of the paper.

\vskip2truecm

\heading{2.  Idea of the construction} \endheading

In this section, we introduce the basic ideas and notation and prove that
the space $Y$ which will be constructed in Sections 3-5 yields solutions
to the problems mentioned in the introduction.

We denote by $c_{00}$ the space of all finitely nonzero (real) sequences.
If $x\in c_{00}$ we denote its co-ordinates by $\{x(j)\}_{j=1}^{\infty}.$
We let $a(x)=\min\{j:  x(j)\neq 0\}$ and $b(x)=\max\{j:x(j)\neq 0\}.$ If
$A$ is a subset of $\bold N$ then $Ax(j)=x(j)\chi_A(j)$ where $\chi_A$ is
the characteristic function of $A.$ If $E_1,E_2$ are subsets of $\bold N$
we write $E_1<E_2$ if $\max E_1<\min E_2.$ We shall also write for
$x,y\in c_{00}$ that $x<y$ if $b(x)<a(y).$ On the hand the natural
co-ordinatewise order on $c_{00}$ will be denoted by $x\le y$, i.e.
$x\le y$ if and only if $x(j)\le y(j)$ for all $j\in\bold N.$ Let
$c_{00}^+=\{x\in c_{00}:x\ge 0\}.$

For $x,y\in c_{00}$ we will write $\langle x,y\rangle=\sum_{j=1}^{\infty}
x(j)y(j).$ We will also the same terminology when $x\in c_{00}^+$ and
$y=\log v$ for some sequence $v\in c_{00}^+$; in this case it will
understood that the pairing can take the value $-\infty$ and that $0\log
0=0.$

By a {\it sequence space} $X$ we will mean a subspace $X$ of the space
$\omega$ of all sequences equipped with a lattice norm $\|\,\|_X$ so
that:\newline (1) $c_{00}\subset X$, \newline (2) If $|x|\le |y|\in X$
then $x\in X$ and $\|x\|_X\le \|y\|_X,$ and \newline (3) If $0\le
x_n\uparrow x$ and $x_n\in X$ with $\sup\|x_n\|_X<\infty$ then $x\in X$
with $\|x\|_X=\sup\|x_n\|_X$ (the Fatou property).  \newline The
canonical basis vectors $\{e_n\}_{n=1}^{\infty}$ then form a
1-unconditional basis for the closure $X_0$ of $c_{00}.$ For convenience
we will write $X^*$ for the K\"othe dual of $X$ which coincides with the
Banach space dual of $X_0.$ We will denote the closed unit ball of a
Banach space $X$ by $B_X.$ We denote the canonical norm on $\ell_p$ by
$\|\,\|_p$ for the cases $p=1$ and $p=\infty.$

Consider a map $\Phi:c_{00} \to \bold R$.  For any $u_1,\ldots,u_n$ we
define
$\Delta_{\Phi}(u_1,\ldots,u_n)=\sum_{i=1}^n\Phi(u_i)-\Phi(\sum_{i=1}^n
u_i).$ $\Phi$ is called {\it quasilinear} if:\newline (4) $\Phi(\alpha
u)=\alpha u$ for $\alpha\in \bold R,$ $u\in c_{00},$ and\newline (5) For
a constant $\delta=\delta(\Phi)$ we have $|\Delta(u,v)|\le
\delta(\|u\|_1+\|v\|_1)$ whenever $u,v\in c_{00}.$

Given a quasilinear map $\Phi$ we can form the twisted sum $Y=\bold
R\oplus_{\Phi} \ell_1$ which is defined to be completion of $\bold
R\oplus c_{00}$ under the quasinorm $$ \|(\alpha,u)\|_{\Phi} =
|\alpha-\Phi(u)| + \|u\|_1.$$ It is readily verified that if $L$ is the
span of the vector $e_0=(1,0)$ then $Y/L$ is isomorphic to $\ell_1.$ This
construction was first used in \cite{13} and \cite{26} where explicit
nontrivial twisted sums of $\bold R$ and $\ell_1$ and hence to deduce
that local convexity is not a three-space property; see also \cite{27}.

\proclaim{Theorem 2.1}Let $\Phi:c_{00}\to\bold R$ be a quasilinear map
and let $Y=\bold R\oplus_{\Phi}\ell_1.$ Then the following conditions are
equivalent:\newline (1) $Y$ contains no basic sequence.\newline (2) If
$Y_0$ is an infinite-dimensional closed subspace of $Y$ then $Y_0$
contains $e_0.$\newline (3) The quotient map $\pi:Y\to \ell_1$ is
strictly singular.\newline (4) $Y$ is topologically minimal.\newline (5)
There is no infinite-dimensional subspace $F$ of $c_{00}$ so that for
some constant $K$ we have $|\Phi(u)|\le K\|u\|_1$ for all $u\in F.$
\newline (6) If $T:\ell_1\to Y$ is a bounded operator then $T$ is
compact.  \newline (7) If $T:Y\to Y$ is a bounded operator then
$T=\lambda I + S$ where $\lambda\in\bold R$ and $S$ is compact.
\endproclaim

\demo{Proof}The equivalence of (1) and (4) is well-known (see Theorem 4.2
of \cite{11} and Theorem 3.2 of \cite{17}, or see \cite{16}).  (2) is
clearly equivalent to (3) and implies (1).  Conversely if (3) fails then
there is an infinite-closed subspace isomorphic to a subspace of
$\ell_1.$ Thus (1)-(4) are all equivalent.

Next we prove (2) implies (5).  Suppose $F$ is an infinite-dimensional
subspace of $c_{00}$ so that $|\Phi(u)|\le K\|u\|_1$ for $u\in E$.  Let
$Y_0$ be the closure of the subspace of all $(0,x)$ for $x\in E.$ Suppose
$(0,x_n)$ converges to $e_0.$ Then $|1-\Phi(x_n)|$ and $\|x_n\|_1$
converge to zero, which is a contradiction.

Next assume (5) and suppose $Y$ contains a basic sequence.  By a
perturbation argument we can suppose it contains a normalized basic
sequence of the form $(\alpha_n,u_n)$ where $u_n\in c_{00}$ By passing to
a subsequence we can suppose that $u_1<u_2<\cdots$ and that $e$ is not in
the closed linear span of $(\alpha_n,u_n).$ It follows that $\pi$ is an
isomorphism on the span of this basic sequence so that for some $K$ we
have:  $$ |\sum_{i=1}^n\alpha_it_i-\Phi(\sum_{i=1}^nt_iu_i)| \le
K\|\sum_{i=1}^nt_iu_i\|_1$$ for all $t_1,\ldots,t_n.$ Let $F_0$ be the
subspace of the linear span of the $(u_n)_{n=1}^{\infty}$ consisting of
all $\sum_{i=1}^nt_iu_i$ where $\sum_{i=1}^n\alpha_it_i=0.$ Then
$|\Phi(u)|\le K\|u\|_1$ for $u\in F_0.$ Thus (5) implies (1).

(3) implies (6).  If $T:\ell_1\to Y$ is bounded then $\pi T$ is strictly
singular and hence compact.  If $(x_n)$ is a sequence in the unit ball of
$\ell_1$ then by passing to a subsequence we can suppose that $\pi Tx_n$
converges.  Hence there exist $y_n\in Y$ so that $(y_n)$ converges and
$\pi Tx_n=\pi y_n$.  But then $Tx_n-y_n\in L$ and so has a convergent
subsequence.

(6) implies (7).  If $T:Y\to Y$ is a bounded operator then since $L$ is
the intersection of the kernels of all continuous linear functionals on
$Y$ we must have $T(L)\subset L.$ Thus $Te=\lambda e$ for some $\lambda.$
Let $S=T-\lambda I;$ then $S= S_0\pi$ where $S_0:Y/L\to Y$ is compact by
(6).

(7) implies (3).  If $\pi$ is not strictly singular, there is a subspace
$Y_0$ of $Y$ of infinite codimension and isomorphic to $\ell_1.$ Hence
there is an isomorphic embedding $V:\ell_1\to Y.$ Then suppose
$V\pi=\lambda I+S$ where $S$ is compact.  Let $\pi_0:Y\to Y/Y_0$ be the
quotient map.  Then $\lambda\pi_0=-S\pi_0$ is compact.  Hence $\lambda=0$
but this contradicts the fact that $V$ is an isomorphism.  \enddemo

\proclaim{Theorem 2.2}If $Y$ satisfies the equivalent conditions of
Theorem 2.1 then any algebraic complement of $L$ satisfies the
Hahn-Banach Extension Property.\endproclaim

\demo{Proof}Let $Z$ be an algebraic complement of $L.$ The continuous
linear functionals on $Z$ separate points so that any linear functional
on a finite-dimensional subspace can be extended continuously to $Z$.
Now let $Z_0$ be a closed infinite-dimensional subspace of $Z$ and
suppose $f$ is a continuous linear functional on $Z_0.$ Let $W$ be the
closure of $Z_0$ in $Y$ and let $f$ denote the extension of $f$ to $W.$
Then $W$ and $f^{-1}(0)$ contain $L$ by (2) and so $f$ factors to a
continuous linear functional on $W/L\subset Y/L$ which is a Banach space.
Hence by the Hahn-Banach theorem $f$ can be extended continuously to $Y$
and hence also to $Z.$\enddemo

\proclaim{Theorem 2.3}If $Y$ satisfies the conditions of Theorem 2.1 then
the topology $\tau$ on $Y$ cannot be the supremum of two vector
topologies $\tau_1,\tau_2$ so that $(Y,\tau_1)$ is nearly convex and
$(Y,\tau_2)$ has trivial dual.\endproclaim

\demo{Proof}Clearly $e_0$ must be in the closure of $\{0\}$ for $\tau_1.$
Let $E$ be the closure of $\{0\}$ for $\tau_2.$ If $e_0\notin E$ then
Theorem 2.1 implies that $E$ is finite-dimensional and that $Y^*$
separates the points of $E$.  Hence $Y=Y_0\oplus E$ for some closed
subspace $Y_0$ of $Y$.  Now $Y_0$ contains no basic sequence and so its
topology is minimal; however $\tau_2$ is Hausdorff on $Y_0$ so that it
must agree with the original topology.  This implies that $Y_0^*=\{0\},$
but in fact $Y_0^*$ is infinite-dimensional.  This contradiction
establishes the theorem.\enddemo

We now review the method of approach to the example.  Theorem 2.1 reduces
the problem to a type of distortion question expressed by (4).  The
recent results of the author \cite{15} show that there is a close
relationship between quasilinear maps on $c_{00}$ and sequence spaces
(see Theorem 6.8 of \cite{15}).  We will explain the connection in the
next section and show how the recent spaces discovered by Gowers and
Maurey (\cite{7} and \cite{9}) enable us to construct a pathological
$\Phi.$

\vskip2truecm \heading{3.  Indicators of sequence spaces}\endheading

We now introduce some ideas from \cite{15}.  Suppose $X$ is a sequence
space.  We define the {\it indicator} $\Phi_X$ (called the {\it entropy
map} in \cite{21}) on $c_{00}$ by $\Phi_X(u)=\langle u,\log x\rangle$
where $u=x^*x$ is the (unique) {\it Lozanovskii factorization} of $u$
i.e.  $x\in B_X^+$ and $x^*\in X^*$ satisfy $\langle x,x^*\rangle
=\|x^*\|_{X^*}=\|u\|_1$ and supp $x$, supp $x^* \subset$ supp $u.$ The
Lozanovskii factorization originates in \cite{19}.

Clearly $\Phi_X(\alpha u)=\alpha\Phi_X(u)$ for $u\in c_{00}.$ Furthermore
if $u,v\in c_{00}$ we also have $$|\Delta (u,v)| \le
\frac4e(\|u\|_1+\|v\|_1)\tag1$$ where $\Delta=\Delta_{\Phi_X}$ (see Lemma
5.6 of \cite{15}).  If $u\in c_{00}^+$ then we can characterize the
Lozanovskii factorization as the solution of an optimization problem so
that $$ \Phi_X(u) =\max_{x\in B_X^+} \langle u,\log x\rangle.\tag2$$ This
idea originates with Gillespie \cite{6}.  Furthermore for
$u_1,\ldots,u_n\in c_{00}^+$ we have the inequalities $$ 0\le
\Delta(u_1,\ldots,u_n) \le
\sum_{i=1}^n\|u_i\|_1\log\frac{S}{\|u_i\|_1}\tag3$$ where
$S=\sum_{i=1}^n\|u_i\|_1$; see \cite{15} Lemma 5.5.

Suppose $f:[1,\infty)\to [1,\infty)$ is any increasing map with $f(1)=1$
and so that $f(t)\le t$ for all $t\ge 1.$ We will say that a sequence
space $X$ has a {\it lower $f$-estimate on blocks} if, whenever
$x_1<x_2<\cdots<x_n\in c_{00}$ then $$ \|x_1+\cdots+x_n\|_X \ge
\frac{1}{f(n)}\sum_{i=1}^n \|x_i\|_X$$ and an {\it upper $f$-estimate on
blocks} if, whenever $x_1<x_2<\cdots<x_n\in c_{00}$ then $$
\|x_1+\cdots+x_n\|_X \le f(n)\max_{1\le i\le n}\|x_i\|_X.$$

\proclaim{Lemma 3.1}Suppose $X$ satisfies an upper $f$-estimate for
blocks.  Then for $u_1<u_2<\cdots<u_n$ in $c_{00}^+$ we have
$$\Delta(u_1,\ldots,u_n) \le \log f(n)(\|u_1\|_1+\cdots +\|u_n\|_1).$$
\endproclaim

\demo{Proof}Let $u_i=x_ix_i^*$ be the Lozanovskii factorizations.  Then
since $f(n)^{-1}(x_1+\cdots+x_n)\in B_X$ we have by (2),
$$\Phi_X(u_1+\cdots+u_n) \ge \langle
\sum_{i=1}^nu_i,\log(f(n)^{-1}\sum_{i=1}^nx_i)\rangle$$ so that the lemma
follows.\enddemo

The following is a special case of Lemma 5.8 of \cite{15}.  Unfortunately
as the referee has pointed out, Lemma 5.8 in \cite{15} is misstated with
the inequality reversed, and in the proof the maximum should be replaced
by the minimum.  This lemma is used in Theorem 5.7 of \cite{15} which is
correct although an inequality is again reversed.  In view of this we
will sketch a simple direct proof.

\proclaim{Lemma 3.2}Suppose $s_1,\ldots,s_n,t_1,\ldots,t_n\ge 0$ and let
$\sum_{i=1}^ns_i=S$ and $\sum_{i=1}^nt_i=T.$ $$ \sum_{i=1}^n
\left(s_i\log\frac{s_i+t_i}{s_i}+t_i\log\frac{s_i+t_i}{t_i}\right)\le
S\log\frac{S+T}{S}+T\log \frac{S+T}{T}.$$ \endproclaim

\demo{Remark}The summand is zero if either $s_i$ or $t_i$
vanishes.\enddemo

\demo{Proof}We will seek to maximize the function
$$u(s_1,\ldots,s_n,t_1,\ldots,t_n) = \sum_{i=1}^n
\left(s_i\log\frac{s_i+t_i}{s_i}+t_i\log\frac{s_i+t_i}{t_i}\right)$$
subject to the constraints $\sum_{i=1}^ns_i=S$ and $\sum_{i=1}^nt_i=T$
and $s_i\ge 0,t_i\ge 0$ for $1\le i\le n.$ By continuity, there is a
point where the maximum is attained.  We can suppose $s_it_i>0$ for $1\le
i\le m$ and $s_it_i=0$ if $m+1\le i\le n$.  By the method of Lagrange
multipliers it is easy to show that $s_i/t_i$ is constant for $1\le i\le
m.$ But then $$ u(s_1,\ldots,s_n,t_1,\ldots,t_n) =
S_0\log\frac{S_0+T_0}{S_0}+T_0\log\frac{S_0+T_0}{T_0}$$ where
$S_0=\sum_{i=1}^ms_i\le S$ and $T_0=\sum_{i=1}^mt_i \le T.$ This
expression is monotone increasing in $S_0$ and $T_0$ and so the result
follows.\enddemo

Let $D=B_{\ell_1}\cap c_{00}^+.$

\proclaim{Lemma 3.3}Suppose $X$ satisfies an upper $f$-estimate on blocks
and suppose $u\in D.$ Let $u=\sum_{i=1}^nu_i$ where $u_1<u_2<\cdots<u_n.$
Let $A$ be any subset of $\bold N$ and let $t=\|Au\|_1.$ Then $$
\Delta(u_1,\ldots,u_n)-(1-t)\log f(n)-\varphi(t) \le
\Delta(Au_1,\ldots,Au_n) \le \Delta(u_1,\ldots,u_n)+\varphi(t),$$ where
$\varphi(t)=t\log\frac1t + (1-t)\log\frac{1}{1-t}(\le \log 2).$
\endproclaim

\demo{Proof}Let $\bold N\setminus A=B.$ Then $$
\Delta(Au_1,\ldots,Au_n,Bu_1,\ldots,Bu_n)=\Delta(Au_1,\ldots,Au_n)+\Delta(Bu_1,
\ldots,Bu_n) +\Delta(Au,Bu).$$ Similarly $$
\Delta(Au_1,\ldots,Au_n,Bu_1,\ldots,Bu_n)=
\Delta(u_1,\ldots,u_n)+\sum_{i=1}^n\Delta(Au_i,Bu_i).$$

Since $\Delta(Bu_1,\ldots,Bu_n),\Delta(Au,Bu)\ge 0$ we deduce
$$\Delta(Au_1,\ldots, Au_n) \le
\Delta(u_1,\ldots,u_n)+\sum_{i=1}^n\Delta(Au_i,Bu_i).$$ Now we use (3)
and Lemma 3.2.  We have $$ \align \sum_{i=1}^n\Delta(Au_i,Bu_i) &\le
\sum_{i=1}^n\left(\|Au_i\|_1\log\frac{\|u_i\|_1}{\|Au_i\|_1}+
\|Bu_i\|_1\log\frac{\|u_i\|_1}{\|Bu_i\|_1}\right)\\ &\le
t\log\frac{1}{t}+(1-t)\log\frac{1}{1-t}.  \endalign $$

For the former inequality we observe that $\Delta(Bu_1,\ldots, Bu_n)\le
\log f(n) \|Bu\|_1.$ Hence $$ \Delta(Au_1,\ldots,Au_n)\ge
\Delta(u_1,\ldots,u_n)-(1-t)\log f(n) -\Delta(Au,Bu),$$ and the second
inequality follows.\enddemo

\proclaim{Lemma 3.4} Suppose $u\in c_{00}^+$ with $\|u\|_1\le 1.$ Suppose
$u=xx^*$ where $x\in B_X^+,x^*\in B_{X^*}^+.$ Then $\Phi_X(u)-\langle
u,\log x\rangle\le \|u\|_1\log\frac1{\|u\|_1}(\le \frac1e).$\endproclaim

\demo{Proof} We can suppose that the supports of $x,x^*$ coincide with
the support of $u.$ Define $Z$ by
$\|z\|_Z=\max(\|z\|_X,\|u\|_1^{-1}\langle |z|,x^*\rangle).$ Then
$\|z\|_X\le \|z\|_Z \le \|u\|^{-1}\|z\|_X$ so that
$\Phi_X(v)+\|v\|_1\log\|u\|_1 \le \Phi_Z(v)\le \Phi_X(v)$ for $v\in
c_{00}^+.$ However $\|x\|_Z\le 1$ and $\|x^*\|_{Z^*}\le \|u\|_1$ so that
$u=xx^*$ is the Lozanovskii factorization for $u.$ Thus
$\Phi_Z(u)=\langle u,\log x\rangle$ and the lemma follows.\enddemo

The next lemma is essentially due to Odell and Schlumprecht \cite{21}.

\proclaim{Lemma 3.5} Given $\epsilon>0$ and $n\in\bold N$ there exists
$\eta>0$ so that if $u_1<u_2<\cdots<u_n$ are in $D$,
$u=\frac1n(u_1+\cdots+u_n)$ are such that
$\delta=\frac1n\Delta(u_1,\ldots,u_n)<\eta$ then for the Lozanovskii
factorizations $u=xx^*$ and $u_i=x_ix_i^*$ we have $\|Au\|_1<\epsilon$
where $A=\{j:y(j)>(1+\epsilon)x(j)\}$ and $y=x_1+\cdots+x_n.$
\endproclaim

\demo{Proof} By Proposition 2.3 of Odell and Schlumprecht, \cite {21},
given $\epsilon>0$ there exists $\nu>0$ so that if $v\in D$ and $z\in
B_X^+$ are such that $\langle v,\log z\rangle > \Phi_X(v)-\nu$ then if
$v=z_0z_0^*$ is the Lozanovskii factorization then $\|Bv\|_1<\epsilon$
where $B=\{j:z_0(j)>(1+\epsilon)z(j)\}.$ Let $\eta=\nu/n.$ Then if
$\delta<\eta$ we have $$ \sum_{i=1}^n(\Phi_X(u_i)-\langle u_i,\log
x\rangle) < \nu$$ and since each term is positive we conclude that
$\|A_iu_i\|_1<\epsilon$ where $A_i=\{j:  x_i(j)>(1+\epsilon)x(j)\}$.
This quickly implies that $\|Au\|_1<\epsilon.$\enddemo

\vskip2truecm

\heading{4.  The Gowers-Maurey space}\endheading

At this point we let $f(x)=\log_2(x+1)$ and introduce as in \cite{9} the
class $\Cal F$ of functions $g:[1,\infty)\to [1,\infty)$ satisfying the
properties:\newline (1) $g(1)=1$ and $g(x)<x$ for $x>1.$ \newline (2) $g$
is strictly increasing and unbounded.  \newline (3)
$\lim_{x\to\infty}x^{-q}g(x)=0$ for any $q>0.$ \newline (4) $x/g(x)$ is
concave and nondecreasing.\newline (5) $g$ is submultiplicative i.e.
$g(xy)\le g(x)g(y)$ for $x,y\ge 1.$ \newline Clearly $f\in \Cal F$ and so
is $\sqrt f.$

Now suppose $X$ is a sequence space.  If $n\in\bold N$ and $\kappa>1$ we
define $\lambda_X(n,\kappa)$ to be the set of $x\in c_{00}^+$ so that
$\|x\|_X=1$ and $x=\frac1n(x_1+\cdots+x_n)$ where $x_1<x_2<\cdots<x_n$
and $\|x_i\|_X\le \kappa$ for $1\le i\le n.$ (Thus $x$ is an
$\ell^n_{1+}$ average with constant $\kappa,$ in the sense of \cite{9}:
note that we restrict to non-negative sequences and to spaces $X$ for
which the canonical basis is unconditional.)

We then define $RIS_X(n;\kappa)$ to be the collection of sequences
$x_1<x_2<\cdots<x_n$ in $c_{00}^+$ satisfying
$x_i\in\lambda_X(M_i,\kappa)$ where $M_1\ge
4\kappa\rho^{-1}2^{36n^2\rho^{-2}}$ and $M_{k+1} \ge
2^{4b(x_k)^2\rho^{-2}}$ for $k\ge 1$ where $\rho=\min(\kappa-1,1).$ We
then define $\Lambda_X(n;\kappa)$ to be the collection of $x\in c_{00}^+$
of the form $x=\|x_1+\cdots+x_n\|_X^{-1}(x_1+\cdots+x_n)$ where
$(x_1,\ldots,x_n)\in RIS_X(n, \kappa).$ This definition differs slightly
but inessentially from that of \cite{9}.  In fact we will only really
require the case $\kappa\ge 2$ when $\rho=1;$ this is in contrast to
\cite{9} where values of $\kappa$ close to one are important.

At the same time if $g\in\Cal F$ we define $\Cal H_X(g;m)$ to be the
collection of $(m,g)$-forms i.e.  $x^*\in\Cal H_X(g;m)$ if and only if
$x^*=g(m)^{-1}(x_1^*+\cdots+x_m^*)$ where $x_1^*<x_2^*<\cdots<x_m^*$ are
in $c_{00}^+$ and $\|x_i^*\|_{X^*}\le 1$ for $1\le i\le m.$

We will require certain lemmas from \cite{9}.

\proclaim{Lemma 4.1}(Lemma 4 of \cite{9}) Suppose
$x\in\lambda_X(N,\kappa)$ and $x^*\in \Cal H_X(g;M)$ where $g\in\Cal F.$
Then $\langle x,x^*\rangle \le \kappa(1+2M/N)g(M)^{-1}.$ \endproclaim

\proclaim{Lemma 4.2}(Lemma 5 of \cite{9}) Suppose $X$ satisfies a lower
$f$-estimate on blocks and $g\in\Cal F$ with $g\ge f^{1/2}.$ Suppose
$N\in\bold N$ and $\kappa>1.$ Suppose $M\ge 2^{36N^2\rho^{-2}}$ and that
$x\in\Lambda(N,\kappa),\ x^*\in\Cal H_X(g,M).$ Then $\langle x,x^*\rangle
\le (\kappa+\rho)f(N)/N\le (\kappa+1)f(N)/N.$\endproclaim

\demo{Remark}For our statement of Lemma 4.2, observe that since $X$ has a
lower $f$-estimate, then for any $\{x_i\}_{i=1}^N\in RIS_X(N,\kappa)$ we
have $\|\sum_{i=1}^Nx_i\|_X \ge N f(N)^{-1}.$\enddemo

Our next lemma is a slight modification of Lemma 7 of \cite{9}.

\proclaim{Lemma 4.3} Suppose $X$ satisfies a lower $f$-estimate on blocks
and $g\in\Cal F$ with $g\ge f^{1/2}.$ Suppose $\kappa\ge 2$ and
$(x_1,\ldots,x_N)\in RIS_X(N,\kappa).$ Let $x=\sum_{i=1}^Nx_i$ and
suppose that for every interval $E$ with $\|Ex\|_X\ge 1$ we have $$
\|Ex\|_X \le \sup\{ \langle Ex,x^*\rangle:\ x^*\in \Cal H_X(g;M),\ M\ge
2\}.\tag *$$ Then $\|x\|_X \le (\kappa+1)N/g(N).$\endproclaim

\demo{Proof}We introduce the length of an interval $E$ as in \cite{9}.
Let $x_i\in\lambda_X(n_i,\kappa)$ for $1\le i\le N.$ Suppose $x_i$ is
written as $\frac1{n_i}\sum_{j=1}^{n_i}x_{ij}$ where
$x_{i1}<x_{i2}<\cdots<x_{in_i}$ and $\|x_{ij}\|_X \le \kappa n_i^{-1}.$
If $E$ is any interval which intersects the support of $\sum_{i=1}^Nx_i$
we let $k\le l$ be the least and greatest indices $i$ such that $Ex_i\neq
0.$ Then we let $p$ be the least index such that $Ex_{kp}\neq 0$ and $q$
the greatest index such that $Ex_{lq}\neq 0. $ Define $\ell(E)= l-k
+qn_l^{-1}-pn_k^{-1}.$ If $E$ does not meet the support of
$\sum_{i=1}^Nx_i$ then $\ell (E)=0.$

Now our hypotheses differ from Lemma 7 of \cite{9} in that we assume (*)
whenever $\|Ex\|_X\ge 1$ while \cite{9} assumes (*) whenever $\ell(E)\ge
1;$ we, however, assume $\kappa\ge 2.$ Our hypotheses imply that (*)
holds if $\ell(E)\ge 2$ since then there exists a least one $x_i$ has
support contained entirely in $E$.  As in \cite{9} let $G(t)=t/g(t)$ for
$t\ge 1$ and $G(t)=t$ for $t\le 1.$ Then if $\kappa n_1^{-1} \le \ell
(E)\le 1$ we have $\|Ex\|_X \le (\kappa+1)G(\ell(E))$ as in \cite{9}.  We
claim the same inequality if $1\le \ell(E)\le 2;$ in fact in this
situation we can see that $E$ intersects the supports of at most three
$x_i$ and so $\|Ex\|_X \le 3\le (\kappa+1)G(\ell(E))$.  The proof can now
be completed by applying Lemma 7 of \cite{9}.  \enddemo

We will now define a Gowers-Maurey space $Z,$ very similar to the
construction in \cite{9}; in fact, essentially the same space is
considered by Gowers in \cite {7} as a counterexample to the hyperplane
problem, and also as a space in which all operators are strictly singular
perturbations of a diagonal map.  We will suppose that $\Cal
P=\{p_k\}_{k=1}^{\infty}$ is an increasing sequence of natural numbers
satisfying $f(p_1)>256,$ $\log\log\log p_k>4p_{k-1}^2,$
$p_k>k^62^{100k^2},$ for all $k.$ We shall also require that
$f(p_{2k})p_{2k}^{-1}\le \frac12 k^{-3}$ which doubtless follows from our
other hypotheses.  For convenience we suppose each $p_k$ is a square.  We
partition $\Cal P=\Cal P_1\cup \Cal P_2$ where $\Cal
P_1=\{p_{2k-1}\}_{k=1}^{\infty}$ and $\Cal
P_2=\{p_{2k}\}_{k=1}^{\infty}.$

Let $\bold Q_+$ denote the countable collection of $u\in c_{00}^+$ which
have only rational coefficients and let $\sigma$ be an injection from the
collection of all finite subsets of $\bold Q_+,$ $\{z_1,z_2,\ldots,z_s\}$
where $z_1<z_2<\cdots<z_s$ to $\Cal P_2$ which satisfies the condition
$\sigma(z_1,\ldots,z_s) \ge 2^{10 b(z_s)^2}.$

We then define $Z$ implicitly by the formula $$ \|x\|_Z =
\max(\|x\|_{\infty},\|x\|_{\alpha},\|x\|_{\beta})$$ where $$
\|x\|_{\alpha} =\sup\{\langle |x|,x^*\rangle:  \ x^*\in\Cal H_Z(f;M),\
M\ge 2\}$$ and $$ \|x\|_{\beta}= \sup\left\{
f(k)^{-1/2}\sum_{i=1}^k\langle |x|,x_i^*\rangle\right\} $$ where the
supremum is over all $k\in \Cal P_1$ and {\it special sequences}
$(x_1^*,\ldots,x_k^*)$ i.e. such that $x_1^*<x_2^*<\cdots<x_k^*,$ with
$x_1^*\in \bold Q_+\cap \Cal H_Z(f;p_{2k})$ and then for $j\ge 1,$
$x_{j+1}^*\in\bold Q_+\cap \Cal H_Z(f;\sigma(x_1^*,x_2^*,\ldots,x_j^*)).$

This implicit definition can be justified by an inductive construction as
in \cite {9}.  Precisely we set $\|x\|_{Z_0}=\|x\|_{\infty}$ for $x\in
c_{00}$ and then define for $N\ge 1,$ $$ \|x\|_{Z_N}
=\max(\|x\|_{Z_{N-1}},\|x\|_{\alpha_{N-1}},\|x\|_{\beta_{N-1}})$$ where
$$ \|x\|_{\alpha_N} =\sup\{\langle |x|,x^*\rangle:  \ x^*\in\Cal
H_{Z_N}(f;M),\ M\ge 2\}$$ and $$ \|x\|_{\beta_N}= \sup\left\{
f(k)^{-1/2}\sum_{i=1}^k\langle |x|,x_i^*\rangle\right\} $$ where the
supremum is over all $k\in \Cal P_1$ and $(x_1^*,\ldots,x_k^*)$ i.e. such
that $x_1^*<x_2^*<\cdots<x_k^*,$ with $x_1^*\in \bold Q_+\cap \Cal
H_{Z_N}(f;p_{2k})$ and then for $j\ge 1,$ $x_{j+1}^*\in\bold Q_+\cap \Cal
H_{Z_N}(f;\sigma(x_1^*,x_2^*,\ldots,x_j^*)).$ It is then easily verified
that $\|\,\|_{Z_N}$ is an increasing sequence of norms, bounded above by
the $\ell_1-$norm and that the sets $H_{Z_N}(f;M)$ also increase in $N.$
We set $\|x\|_Z=\lim_{N\to\infty}\|x\|_{Z_N}.$

We emphasize that this space is an unconditional version of the
counterexample constructed in \cite{9}, but shares some of the same
features.  We will need versions for $Z$, of certain lemmas proved in
\cite{9} for the Gowers-Maurey space.  Fortunately the same basic
techniques go through more or less unchanged.

Let us note first that $Z$ satisfies a lower $f$-estimate.  This follows
immediately from the definition of $\|x\|_{\alpha}.$ We also note that,
by induction, it follows that $\|e_n\|_Z=1$ for all $n.$

\proclaim{Lemma 4.4}Suppose $(x_j)_{j=1}^n\in RIS_Z(n;\kappa)$ where
$\kappa\ge 1.$ Then $\|\sum_{j=1}^nx_j\|_{\infty}<1.$ \endproclaim

\demo{Proof}We have $x_j\in\lambda_Z(M_j,\kappa)$ where $M_j\ge 4\kappa$
by the definition of $RIS_Z(n;\kappa).$ Hence $\|x_j\|_{\infty}\le
M_j^{-1}\kappa<1$ and the lemma follows.\enddemo

It now follows as in Lemma 10 of \cite{9}:

\proclaim{Lemma 4.5} Suppose $\kappa\ge 2.$ Suppose $N\in \Cal P_2$ and
$\log N\le n\le \exp N.$ Then if $\{x_1,\ldots,x_n\}\in RIS_Z(n,\kappa)$
we have $\|\sum_{i=1}^nx_i\|_Z \le (\kappa +1)nf(n)^{-1}.$\endproclaim

\demo{Proof}The key point proved in \cite{9}, Lemma 9, is that there
exists $g\in\Cal F$ with $f^{1/2}\le g\le f$ such that $g(x)=f(x)$ for
$\log N\le x\le \exp N$ and $g(k)=f^{1/2}(k)$ when $k\in \Cal P_1.$ Thus
if $x\in c_{00}$ and $\|x\|_Z>\|x\|_{\infty}$ then $$\|x\|_Z =
\sup\{\langle Ex,x^*\rangle; \ x^*\in\Cal H_Z(g;M), \ M\ge 2\rangle.$$

Now, by the preceding lemma if $x=\sum_{j=1}^nx_j$ and $E$ is any
interval then $\|Ex\|_{\infty}<1.$ We can therefore apply Lemma 4.3 to
obtain the result.\enddemo

The next lemma is simply a cruder form of Lemma 11 from \cite{9}.

\proclaim{Lemma 4.6} Suppose $\kappa\ge 2$ and $N\in \Cal P_2.$ If $x\in
\Lambda_Z(N,\kappa)$ then $x \in \lambda(\sqrt N,2(\kappa+1)).$
\endproclaim

\demo{Proof}Suppose $\{x_i\}_{i=1}^N\in RIS_Z(N,\kappa)$ and that
$x=\|\sum_{i=1}^Nx_i\|_Z^{-1}\sum_{i=1}^Nx_i.$ We break $[1,N]$ into
$\sqrt N$ intervals $E_j$ each of length $\sqrt N,$ which is an integer
by hypothesis.  Note that $\{x_i\}_{i\in E_j}\in RIS_Z(\sqrt N,\kappa).$
If $y_j=\sum_{i\in E_j}x_i$ then, by Lemma 4.5, $\|y_j\|_Z \le
(\kappa+1)\sqrt N.$ Also $\|\sum_{j=1}^{N}x_j\|_Z \ge N/f( N ),$ by the
lower $f$-estimate on $X.$ Now $x=\frac1{\sqrt N}(\sum_{j=1}^{\sqrt
N}z_j)$ where $ z_j=(\|\sum_{i=1}^Nx_i\|_Z)^{-1}\sqrt N y_j$.  But
$\|z_j\|_Z \le (\kappa+1)(Nf(N))/(Nf(\sqrt N))\le 2(\kappa+1).$ \enddemo

Our next result is a modification of Lemma 12 of \cite{9}.  In fact, this
Lemma appears to be incorrectly stated in \cite{9} and so some
modification is necessary.  In the proof of the lemma in \cite{9} it is
claimed without justification that $\{x_1,\ldots,x_k\}$ is a RIS of
length $k$ and constant $1+\epsilon.$ For the applications some
modification similar to that given below seems adequate, however.

\proclaim{Lemma 4.7} Let us suppose $\kappa\ge 2.$ Suppose $k\in \Cal
P_1$ with $f(k)>100\kappa^2.$ Suppose $E_1,\ldots,E_k$ are intervals with
$E_1<E_2<\cdots<E_k.$ Let $\{x_1^*,\ldots,x_k^*\}$ be a special sequence
with supp $x_j^*\subset E_j.$ Let $M_1=p_{2k}$ and
$M_{j+1}=\sigma(x_1^*,\ldots,x_j^*)$ for $1\le j\le k-1.$ Let $A$ be any
subset of $\{1,2,\ldots,k\}$ and suppose for each $j\in A$ we have
$x_j\in c_{00}^+$ with supp $x_j\subset E_j,$ so that $x_j,x_j^*$ are
disjoint and $x_j\in\Lambda(M_j,\kappa).$ Then $$ \|\sum_{i\in A}x_i\|_Z
\le 16\kappa k f(k)^{-1}.$$ \endproclaim

\demo{Proof}We have $x_j\in\lambda_Z(\sqrt {M_j}, 4\kappa),$ by Lemma
4.6.  Note that $\sqrt {M_1}=\sqrt {p_{2k}}\ge 4\kappa 2^{36k^2}.$ We
also have $\sqrt {M_{j+1}} > 2^{4b(x_j^*)^2}.$

Now assume $A$ contains no two consecutive integers.  Then if $j\in A$ we
have $\sqrt {M_j}\ge 2^{4b(x_{j-2})^2}$ for $j\ge 2$ and so
$\{x_j\}_{j\in A}\in RIS_Z(|A|,4\kappa)$.  As in \cite{9} we use Lemma
4.3.

Note first that there exists $h\in\Cal F$ with $\sqrt f\le h\le f,$ so
that $h(n)=\sqrt {f(n)}$ if $n\in \Cal P_1\setminus\{k\}$ while
$h(n)=f(n)$ if $n\in \Cal P_2\cup\{k\}.$ This fact follows from Lemma 9
of \cite{9}.

Let $x=\sum_{i\in A}x_i$ and suppose, for some interval $E$ we have
$\|Ex\|_Z\ge 1,$ and $$ \|Ex\|_Z >\sup\left\{\langle Ex,x^*\rangle:\
x^*\in\Cal H_Z(h,m),\ m\ge 2\right\}.$$ Since $h\le f$ this implies that
$\|Ex\|_Z >\|Ex\|_{\alpha}.$ On the other hand, since $\{x_j\}_{j\in
A}\in RIS_Z(|A|,4\kappa)$ we can apply Lemma 4.4 to deduce that
$\|Ex\|_Z>\|Ex\|_{\infty}.$ The conclusion is that
$\|Ex\|_Z=\|Ex\|_{\beta}.$ Thus there is a special sequence
$\{z_1^*,z_2^*,\ldots,z^*_l\},$ with $l\in\Cal P_1,$ so that $$ \|Ex\|_Z
=f(l)^{-1/2} \langle Ex, \sum_{i=1}^lz_i^*\rangle.$$ However,
$f(l)^{1/2}=h(l)$ unless $l=k.$ We conclude $l=k$ and $$ 1\le \|Ex\|_Z
\le f(k)^{-1/2}\sum_{i\in A}\sum_{j=1}^k\langle x_i,z_j^*\rangle.$$

Let $t$ be the greatest integer so that $z_t^*=x_t^*$ (with $t=0$ if no
such integer exists).  If $i<t$ it is clear that $\langle
x_i,z_j^*\rangle =0$ for all $j.$ Similarly if $j\le t$ it is also clear
that $\langle x_i,z_j^*\rangle=0$ for all $i.$ If $i=t,$ then $\langle
x_i,z_j^*\rangle =0$ unless $j=t+1$ when of course $\langle
x_t,z_{t+1}^*\rangle \le 1.$ If $t+1\le i\in A$ and $t+1\le j\le k$ then,
unless $t+1=i=j$ we have $x_i\in \Lambda_Z(M_i,\kappa)$ and $z_j^*\in\Cal
H_Z(g;M'_j)$ where $M_i,M'_j\in\Cal P_2$ are not equal.  It follows from
the separation conditions on $\Cal P_2$ that we can apply either Lemma
4.1 or Lemma 4.2; if $M'_j<M_i,$ then by Lemma 4.1, $$ \langle
x_i,z_j^*\rangle \le 24\kappa f(M'_j)^{-1} \le 24\kappa f(p_{2k})^{-1},
$$ or if $M'_j>M_i$, then $M'_j\ge 2^{36M_i^2}$ and by Lemma 4.2, $$
\langle x_i,z_j^*\rangle \le 2\kappa f(M_i)/M_i\le 2\kappa
f(p_{2k})p_{2k}^{-1}.$$ In either case we have $\langle x_i,z_j^*\rangle
\le \kappa k^{-2} $. If $i=j=t+1$ then $\langle x_i,z_j^*\rangle\le 1.$

Hence $$\langle \sum_{i\in A}x_i,\sum_{j=1}^kz_j^*\rangle \le 2+\kappa\le
3\kappa.$$

This implies that $$\|Ex\|_Z \le 3\kappa f(k)^{-1/2}<\frac3{10}$$
contrary to assumption.  The conclusion from Lemma 4.3 is then that
$$\|x\|_Z \le 8\kappa |A| h(|A|)^{-1} \le 8\kappa k f(k)^{-1}.$$ The
general result follows by splitting $A$ into two subsets obeying the
condition that no two consecutive integers are contained in either.
\enddemo

\vskip2truecm

\heading{5.  The main result }\endheading

We now let $X=Z^*$ and consider the indicator $\Phi_X$.  We will need the
elementary fact, which follows from duality, that $X$ satisfies an upper
$f$-estimate, i.e. if $x_1<x_2<\cdots<x_n\in c_{00}$ then
$\|x_1+\cdots+x_n\|_X \le f(n)\max_{1\le i \le n}\|x_i\|_X.$ It also
follows from the definition of $Z$ that if $x_1,\ldots,x_n$ is a special
sequence (with $n\in\Cal P_1$) then $\|x_1+\cdots+x_n\|_X \le
f(n)^{1/2}.$

Our main result, which combined with the results of Section 2 establishes
Theorems 1.1, 1.3 and 1.4, is the following:

\proclaim{Theorem 5.1}For every infinite-dimensional subspace $G$ of
$c_{00}$ we have $\sup\{|\Phi_X(u)|:\ \|u\|_1=1,\ u\in
G\}=\infty.$\endproclaim

\demo{Remark}The following proof has been substantially simplified
according to a suggestion of B. Maurey.\enddemo

\demo{Proof}We will start from the assumption that there is a subspace
$G$ of infinite dimension so that $|\Phi_X(u)|\le K\|u\|_1$ for $u\in G$.
We may suppose that if $u\in G$ then $\langle u,\chi\rangle=0$ where
$\chi$ is the constantly one sequence.  Then by induction we can pick
$\xi_1<\xi_2<\xi_3<\cdots$ in $G$ with $\|\xi_j\|_1=2.$ We split $\xi_i$
into positive and negative parts $\xi_i=\xi'_i-\xi''_i,$ where
$\xi'_i,\xi''_i$ are disjoint and nonnegative.  Then $\xi'_i,\xi''_i\in
D.$ We let $R$ be the union of the supports of the $\xi'_i$ and $S$ be
the union of the supports of the $\xi''_i.$ Let $W$ be the linear span of
$\{|\xi_i|\}_{i=1}^{\infty}.$

Notice first that $X$ satisfies an upper $f$-estimate on blocks where
$f(x)=\log_2(x+1).$ If $\gamma>0$ and $n\in\bold N$ we define
$\Gamma(n,\gamma)$ to be the set of $w\in D$ such that there exist
$w_1<w_2<\cdots<w_n\in D$ with $w=\frac1n(w_1+\cdots+w_n)$ and
$\frac1n\Delta(w_1,\cdots,w_n)< \gamma.$

\proclaim{Lemma 5.2}Given any $m,n\in\bold N$ and $\delta>0$ there exists
$w\in W\cap \Gamma(n,\delta)$ with $m< a(w).$ \endproclaim

\demo{Proof}For $n\in\bold N$ let $c_n$ be the infimum of all constants
$\gamma$ so that if $m\in\bold N$ there exists $w\in W\cap
\Gamma(n,\gamma)$ with $m<a(w).$ It is easy to see that $ c_{np}\ge
c_{n}+c_{p}$ for any $n,p$ and that from Lemma 3.1 $c_n\le \log f(n).$
Hence $pc_n\le c_{n^p} \le \log f(n^p)$ and so letting $p\to\infty$ we
obtain $c_n=0$ for all $n$ and the lemma follows.\enddemo

We now turn to estimates on the Lozanovskii factorization of $w\in
\Gamma(n,\delta).$

\proclaim{Lemma 5.3} For fixed $n$ and $0<\epsilon<\frac12$ there exists
$\eta>0$ so that if $w\in\Gamma(n,\eta)$ and $w=xx^*$ is the Lozanovskii
factorization of $w,$ then there exists $A\subset [a(w),b(w)]$ with
$\|Aw\|_1>1-\epsilon$ and such that $Ax^*/\|Ax^*\|_Z\in \lambda_Z(n,2).$
\endproclaim

\demo{Proof}If $w\in\Gamma(n,\delta)$ then $w=\frac1n\sum_{i=1}^nw_i$
where $w_1<w_2<\cdots<w_n\in D$ are such that
$\frac1n\Delta(w_1,\ldots,w_n)\le \delta.$ Let $w_i=x_ix_i^*$ be the
Lozanovskii factorizations of each.  Let $y=x_1+\cdots+x_n$.  If $c>1$
let $A=\{j:  y(j)\le cx(j),\ x(j)>0\}.$ Then $Ax^* \le
cn^{-1}A(x^*_1+\cdots+x^*_n)$ and hence if $A_i=A\cap[a(w_i),b(w_i)]$
then $\|A_ix^*\|_Z\le c/n.$ Now $\|Ax^*\|_Z\ge \|Aw\|_1$ and so
$\|Ax^*\|_Z^{-1}Ax^*\in\lambda_Z(n,c')$ where $c'\le c\|Aw\|_1^{-1}.$
Now, according to Lemma 3.5, if $\delta>0$ is sufficiently small we can
choose $c$ close enough to $1$ so that the conclusions follow.  \enddemo

Using the preceding lemma we describe a construction.  Suppose $N\in\Cal
P_2$ and $\epsilon>0.$ Then given any $m\in\bold N$ and any
$M_1>2^{36N^2+4}$ we can construct two sequences $\{w_j\}_{j=1}^N$ and
$\{\zeta_j\}_{j=1}^N$ and a sequence of integers $(M_j)_{j=1}^N$ so
that:\newline (1) $ m<a(w_1),$ \newline (2) $
w_1<\zeta_1<w_2<\zeta_2<\cdots<w_N<\zeta_N,$\newline (3) $
w_j\in\Gamma(M_j,\eta_j)\cap W$ where $0<\eta_j<\epsilon$ is sufficiently
small so that there exists $A_j\subset [a(w_j),b(w_j)]$ with
$\|A_jw_j\|_1>1-\epsilon$ and $z_j=\|A_jx_j^*\|_Z^{-1}A_jx_j^*\in
\lambda_Z(M_j,2)$ where $w_j=x_jx_j^*$ is the Lozanovskii factorization
of $w_j.$\newline (4) $\zeta_j \in\lambda_Z(M_j,2)$ \newline (5)
$M_{j+1}>2^{4b(\zeta_j)^2}.$

We will call the resulting sequence $\{w_j\}_{j=1}^N$ an
$(N,\epsilon)-${\it sequence} and $w=\frac1N(w_1+\cdots+w_N),$ the
associated $(N,\epsilon)-${\it average.} The sequence
$\{\zeta_j\}_{j=1}^N$ is called the {\it ballast sequence}; it is present
is simply for technical reasons to provide ballast in the argument.  Let
$H$ be the union of the supports of the ballast sequence.

\proclaim{Lemma 5.4} Suppose $\{w_1,\ldots,w_N\}$ is an
$(N,\epsilon)-$sequence as above with associated $(N,\epsilon)-$average
$w$ and ballast $\{\zeta_j\}_{j=1}^N$.  Then there is a subset $A$ of
$[a(w),b(w)]$ and $x\in \Cal H_Z(f;N)\cap\bold Q_+$ with $\supp x\subset
\supp w$ so that:\newline (6) $\|Aw\|_1>1-\epsilon$ \newline (7) If
$B\subset A$ there exists $z\in\Lambda_Z(N,4)$ supported in $B\cup H$ so
that $Bw\le 10xz.$\newline (8) If $B\subset A$ then $\langle Bw,\log
x\rangle>\Phi_X(Bw)-4.$ \endproclaim

\demo{Proof} Notice that $y=\frac1{f(N)}(x_1+\cdots+x_N)\in \Cal
H_Z(f;N)$ and $\|y\|_X\le 1,$ since $X$ has an upper $f$-estimate.
Choose $x$ with rational coefficients so that $\frac12y\le x\le y.$ Let
$A=A_1\cup\cdots\cup A_N$ so that (6) immediately holds.

We recall that $z_j\in\lambda_Z(M_j,2)$ (condition (3)) for $1\le j\le
N.$ It follows easily that if $B$ is a subset of $A$ then we can find
$0\le \alpha_j\le 1$ so that $\|Bx_j^*+\alpha_j\zeta_j\|_Z=1$ and then
$Bx_j^*+\alpha_j\zeta_j \in \lambda_Z(M_j,4).$ The sequence
$\{Bx_j^*+\alpha_j\zeta_j\}_{j=1}^N$ thus belongs to $RIS_Z(N,4)$ (since
$M_1>2^{36N^2+4}$) and so $\|\sum_{j=1}^N(Bx_j^*+\alpha_j\zeta_j)\|_Z \le
5N/f(N),$ from Lemma 4.5.

Let $z$ be the normalized vector
$\beta(\sum_{j=1}^NBx_j^*+\alpha_j\zeta_j)$ where, by the above, $\beta
\ge \frac{f(N)}{5N}.$ Then $z\in\Lambda_Z(N,4)$ and $xz\ge \frac12yz\ge
\frac1{10}Bw.$ This proves (7).

For (8) we notice that Lemma 3.4 now implies that $\Phi_X(Bw) - \langle
Bw,\log x\rangle \le 10/e<4.$

\enddemo

Let us suppose that $n\in \Cal P_1$ is fixed and large, say
$f(n)>\exp(8K+4000),$ and let $\epsilon=(\log f(n))^{-1}.$ Let
$M_1=p_{2n}$; we can construct a $(M_1,\epsilon)$-sequence
$\{w_{1j}\}_{j=1}^{M_1}$ with $(M_1,\epsilon)-$average
$w_1=M_1^{-1}\sum_{j=1}^{M_1}w_{1j},$ and ballast
$\{\zeta_{1j}\}_{j=1}^{M_1}.$ Let $x_1\in\bold Q_+\cap \Cal H_Z(f;M_1)$
and $A_1\subset [a(w_1),b(w_1)]$ be such that the conclusions of Lemma
5.4 hold.

Next let $M_2=\sigma(Rx_1)$ and and construct a $(M_2,\epsilon)-$sequence
$\{w_{2j}\}_{j=1}^{M_2}$ with associated $(M_2,\epsilon)-$average $w_2$
and ballast $\{\zeta_{2j}\}_{j=1}^{M_2}$ so that $\zeta_{1M_1}<u_2.$
Repeating this construction for $n$ steps we obtain sequences
$(w_{ij})_{j=1}^{M_i},(\zeta_{ij})_{j=1}^{M_i}$ for $i=1,2,\ldots,n,$
$(w_i)_{i=1}^n,$ $(M_i)_{i=1}^n,$ $(A_i)_{i=1}^n$ and $(x_i)_{i=1}^n$ so
that:\newline (9) $(w_{ij})_{j=1}^{M_i}$ is an $(M_i,\epsilon)$-sequence
with associated $(M_i,\epsilon)$-average $w_i$ and ballast
$\{\zeta_{ij}\}_{j=1}^{M_i}$ for $1\le i\le n.$ \newline (10)
$w_1<\zeta_{1M_1}<w_2<\zeta_{2M_2}\ldots<w_n<\zeta_{nM_n}.$\newline (11)
$A_i\subset [a(w_i),b(w_i)]$ for $1\le i\le n$ and
$\|A_iw_i\|_1>1-\epsilon_i.$ \newline (12) $\supp x_i\subset \supp w_i,$
$x_i\in \Cal H_Z(f;M_i)\cap\bold Q_+,$ and so $\|x_i\|_X\le 1.$\newline
(13) $\langle Bw_i,\log x_i\rangle>\Phi_X(Bw_i)-4,$ whenever $B\subset
A_i.$\newline (14) For any $B\subset A_i$ there exists
$z\in\Lambda_Z(M_i,4)$ with $Bw_i\le 10x_iz.$\newline (15)
$M_{i+1}=\sigma(Rx_1,\ldots,Rx_i)$ for $1\le i\le n-1.$

We also have:\newline (16) $(Rx_1,\ldots,Rx_n)$ is a special sequence of
length $n$ in $X=Z^*.$

Let $H_i$ be the union of the supports of the ballast at the $i$th. step.
Let $A=\cup_{i=1}^nA_i$ and then set $P=A\cap R$ and $Q=A\cap S.$ We also
define $u_i=2Rw_i,v_i=2Sw_i$ (so that $u_i,v_i\in D$) and then set
$u=\frac1n(u_1+\cdots+u_n),\ v=\frac1n(v_1+\cdots+v_n)$ and
$w=\frac1n(w_1+\cdots+w_n).$

If we set $x=(f(n))^{-1/2}\sum_{i=1}^nRx_i$ then $\|x\|_X\le 1,$ since by
(16) $\{Rx_1,\ldots,Rx_n\}$ is a special sequence.  Hence, using (13)
above, $$ \Phi_X (Pw) \ge \langle Pw,\log x\rangle \ge
\frac1n\sum_{i=1}^n\Phi_X (Pw_i) -\frac12\log f(n)\|Pw\|_1 -4.$$ It
follows that $ \frac1n\Delta(Pu_1,\ldots,Pu_n) \le \frac12\log f(n)+8.$
Now $\frac1n\sum_{i=1}^n\|Pu_i\|_1>1-2\epsilon$ so that by Lemma 3.3, and
the choice of $\epsilon,$ $$ \frac1n\Delta(u_1,\ldots,u_n) \le
\frac12\log f(n) + 11.\tag 17$$

On the other hand, by Lemma 5.4 we can find $z_i\in\Lambda_Z(M_i,4)$
supported on $((\supp w_i)\cap Q)\cup H_i$ so that $Qw_i \le 10x_iz_i.$
At this point we can invoke Lemma 4.7.  Let
$E_i=[a(w_i),b(\zeta_{i,M_i})]$ and notice that $Rx_i,z_i$ are both
supported in $E_i,$ but are disjoint.  Since $f(n)\ge 1600,$
$Rx_1,\ldots,Rx_n$ is a special sequence and $z_i\in\Lambda_Z(M_j,4)$
where $M_1=p_{2n}$ and $M_{j+1}=\sigma(x_1,\ldots,x_j)$ for $1\le j\le
n-1$ we can conclude that $ \|\sum_{j=1}^N z_j\|_Z \le 64nf(n)^{-1}.$ At
the same time, by the upper $f$-estimate on $X,$ $\|\sum_{i=1}^nx_i\|_X
\le f(n).$ Now we have $$ \frac{1}{640}Qw \le
(\frac{1}{f(n)}\sum_{i=1}^nx_i)(\frac{f(n)}{64n}\sum_{i=1}^nz_i)$$ and we
can apply Lemma 3.4 again to deduce that $$ \Phi_X(Qw) \le
\sum_{i=1}^n\langle Qw,\log x_i\rangle -\log f(n)\|Qw\|_1 +640,$$ and so,
since $\|Qw\|_1>\frac12-\epsilon,$ $$ \frac1n\Delta(Qw_1,\ldots,Qw_n) \ge
\log f(n)\|Qw\|_1 -640\ge \frac12\log f(n)-641.$$ Now recall that
$Qv_i=2Qw_i.$ Hence $ \frac1n\Delta(Qv_1,\ldots,Qv_n) \ge \log f(n)-1282$
and we can apply Lemma 3.3 to deduce that $$
\frac1n\Delta(v_1,\ldots,v_n) \ge \log f(n) -1283.\tag 18$$

Notice that $u_i-v_i\in G$ for $1\le i\le n.$ Now we have that
$|\Phi_X(u_i-v_i)|\le 2K$, for $1\le i\le n$ and $|\Phi_X(u-v)|\le 2K.$
Hence $|\Phi_X(u_i)-\Phi_X(v_i)|\le 2K+8e^{-1}\le 2K+3$ for $1\le i\le n$
and similarly $|\Phi_X(u)-\Phi_X(v)|\le 2K+3.$ This implies that $$
\frac1n\Delta(v_1,\ldots,v_n)-\frac1n\Delta(u_1,\ldots,u_n)\le 4K+6.$$
Combining with (17) and (18) gives that $ \log f(n) \le 8K+2600$ which
contradicts our initial choice of $n$ and completes the proof.\enddemo

It is perhaps worth noting at this point that it is very simple to modify
our example so that Theorem 1.1 holds with $L$ of any specified
dimension.

\proclaim{Theorem 5.5}For any $n\in\bold N$ there is a quasi-Banach space
$Y^{(n)}$ with a subspace $L$ of dimension $n$ so that $Y/L$ is
isomorphic to $\ell_1$ and if $Y_0$ is a closed infinite-dimensional
subspace of $Y^{(n)}$ then $L\subset Y.$\endproclaim

\demo{Proof}Let $A_k=\{nj+k\}_{j=0}^{\infty}\subset \bold N,$ for
$k=1,2,\ldots,n.$ Define $S_k:c_{00}\to c_{00}$ by
$S_ku=\sum_{j=0}^{\infty}u(j)e_{nj+k}.$ Define $\Phi:c_{00}\to
\ell_{\infty}^n$ by $\Phi(u)=\{\Phi_X(S_ku)\}_{k=1}^n.$ Then let
$Y^{(n)}$ be the completion of $\ell_{\infty}^n\oplus c_{00}$ under the
quasi-norm $$ \|(\xi,u)\|_{\Phi}=\|\xi-\Phi(u)\|_{\infty}+ \|u\|_1.$$ Let
$L$ be the space of all $(\xi,0)$ for $\xi\in\ell_{\infty}^n.$ Clearly
$Y^{(n)}/L$ is isomorphic to $\ell_1.$ Now suppose $Y_0$ is an
infinite-dimensional subspace so that $Y_0\cap L$ is a proper subspace of
$L$.  Then there is a non-trivial linear functional $f$ on
$\ell_{\infty}^n$ so that $Y_0\cap L\subset Z=f^{-1}(0).$ Suppose
$f(\xi)=\sum_{k=1}^n\beta_k\xi_k.$ It is easy to verify that $Y/Z$ is
isomorphic to the completion of $\bold R\oplus c_{00}$ under the
quasi-norm $\|(\alpha,u)\|_{\Psi}=|\alpha-\Psi(u)|+\|u\|_1$ where
$\Psi(u)=\sum_{k=1}^n\beta_k\Phi_X(S_ku).$ However there is a constant
$K$ depending only on $\beta_1,\ldots,\beta_n$ so that
$|\Psi(u)-\Phi_X(\sum_{k=1}^n\beta_kS_ku)|\le K\|u\|_1.$ It follows
easily that $\Psi$ is unbounded on every infinite-dimensional subspace of
$c_{00}$ and hence that $(Y_0+Z)/Z$ must contain $L/Z$ which is a
contradiction to the fact that $Y_0\cap L$ is contained in $Z.$\enddemo

\vskip2truecm

\heading{6.  Some final remarks}\endheading

In this short final section we will present a proof of Theorem 1.2, which
first appeared in \cite{14}, a reference which may not be readily
available.  Our proof here is slightly shorter.  We begin with a lemma:

\proclaim{Lemma 6.1}Suppose $X$ is a quasi-Banach space with a dense
subspace $V$ with (HBEP).  Suppose $L=\{x\in X:\ x^*(x)=0\ \forall x^*\in
X^*\}.$ Then:\newline (1) If $L=\{0\},$ so that $X$ has a separating dual
then $X$ is locally convex.\newline (2) If $X$ contains a basic sequence
then $X$ is locally convex.  \newline (3) If $M$ is a closed subspace of
$L$ then $X/M$ has a dense subspace with (HBEP).  \endproclaim

\demo{Proof}(1) (cf.\cite{11}) Let $\|\,\|_c$ be the Banach envelope norm
on $X$ i.e.  $\|x\|_c=\sup\{|x^*(x)|:\ \|x^*\|\le 1\}.$ If $X$ is not
locally convex we may choose $v_n\in V$ with $\|v_n\|_c\le 4^{-n}$ and
$\|v_n\|=1.$ Pick any $x\in V$ and consider the sequence
$w_n=v_n+2^{-n}x.$ Then (see Theorem 4.7 of \cite{16}) there is a
subsequence $(w_{p_n})$ which is a Markushevich basis for its closed
linear span in $X$.  Pick $n_0$ large enough so that $x\notin
[w_{p_k}:k\ge n].$ Then by (HBEP) for $V$ there is a linear functional
$x^*\in X^*$ with $x^*(w_{p_k})=0$ for $k\ge n$ but $x^*(x)=1.$ However
$\lim_{n\to\infty}\|x-2^{n}w_n\|_c =0$ so that $x^*(x)=0$ contrary to
hypothesis.

(2) Pick any $u\in L:$ we will show $u=0.$ Assume then that $u\neq 0.$
Suppose $w\in V$ is nonzero, and $u,w$ are linearly independent.  Since
$X$ contains a basic sequence and $V$ is dense in $X$ we can apply
standard perturbation arguments to suppose that we have a bounded basic
sequence $(x_n)$ with $x_n\in n(u+w)+V$, say $x_n=n(u+w)+v_n$ where
$v_n\in V.$ Then there exists $n_0$ so that $[u,w]\cap [x_n]_{n\ge
n_0}=\{0\}.$ Thus there is a bounded linear functional $f$ on the span
$Y$ of $u,w$ and $[x_n]_{n\ge n_0}$ with $f(u)=1,$ $f(w)=0$ and
$f(x_n)=0$ for $n\ge n_0.$ Since $V$ has (HBEP) there is a bounded linear
functional $x^*$ on $X$ with $x^*(v)=f(v)$ for $v\in V\cap Y.$ Thus
$x^*(w)=0$ and $x^*(v_n)=-n;$ also $x^*(u)=0$ since $u\in L.$ Hence
$x^*(x_n)= -n$ contradicting the boundedness of $x^*.$ Now since
$L=\{0\}$ we can apply (1) to deduce that $X$ is locally convex.

(3) Let $\pi:X\to X/M$ be the quotient map; we show $\pi(V)$ has (HBEP).
Indeed if $E\subset \pi(V)$ is a subspace and $f$ is a continuous linear
functional on $E$ then we can find $x^*\in X^*$ so that $x^*(v)=f(\pi v)$
for $v\in \pi^{-1}E\cap V.$ But then $x^*(x)=0$ if $x\in M\subset L$ so
that $x^*$ factors to a linear functional on $X/M.$\enddemo

\proclaim{Theorem 6.2}Suppose $X$ is a decomposable quasi-Banach space.
If $X$ has a dense subspace $V$ with (HBEP) then $X$ is locally
convex.\endproclaim

\demo{Proof}Let $P$ be a bounded projection on $X$ so that both $P$ and
$Q=I-P$ have infinite rank.  If $L$ is defined as in the previous lemma
then $L$ is clearly invariant for $P.$ From the hypotheses, $X^*$ has
infinite dimension and hence so has $X/L.$ Therefore either $P(X)/P(L)$
or $Q(X)/Q(L)$ has infinite dimension.  Suppose the former; then consider
$X/P(L)$ which has a dense subspace with (HBEP) by Lemma 6.1 (3).  Then
$P(X)/P(L)$ is isomorphic to a subspace of $X/L$ which has separating
dual; since it has infinite dimension, it contains a basic sequence.  By
Lemma 6.1 (2) this implies that $X/P(L)$ is locally convex and hence that
$Q(X)$ is locally convex.  But now $X$ itself must contain a basic
sequence and Lemma 6.1 (2) shows that $X$ is locally convex.\enddemo

Let us conclude by mentioning that in \cite{14} we raised the question of
whether every quasi-Banach space $X$ with separating dual has a weakly
closed subspace $W$ and a bounded linear functional $f$ on $W$ which
cannot be extended to $X.$ We proved that this is equivalent to the
following:

\proclaim{Problem}Suppose $X$ is a quasi-Banach space with separating
dual and suppose that every quotient $X/E$ by an infinite-dimensional
subspace $E$ is locally convex.  Is $X$ locally convex?\endproclaim

Of course our main example $Y$ has every quotient $Y/E$ by an
infinite-dimensional subspace locally convex, but fails to have a
separating dual.

\enddocument